\documentclass[reqno, 12pt]{amsart}

 \usepackage{color}

\newtheorem{theorem}{Theorem}[section]

\newtheorem{lemma}[theorem]{Lemma}

\theoremstyle{definition}
\newtheorem{remark}[theorem]{Remark}

\theoremstyle{definition}
\newtheorem{definition}[theorem]{Definition}

\theoremstyle{definition}
\newtheorem{assumption}[theorem]{Assumption}

\newcommand{\mysection}[1]{\section{#1}
\setcounter{equation}{0}}

 \newcommand{\dashint}{\,\,\hbox{\rm\bf--}\kern-.98em\int}
 
\def\bB{\mathbb{B}}
\def\bC{\mathbb{C}}
\def\bM{\mathbb{M}}
\def\bR{\mathbb{R}}
\def\bZ{\mathbb{Z}}

\def\cF{\mathcal{F}}

\def\cM{\mathcal{M}}

\def\cL{\mathcal{L}}

\title[Elliptic equations] 
{Second-order elliptic  equations with
variably partially  VMO coefficients}
\author{N.V. Krylov}
\thanks{The work   was partially supported by
NSF Grant DMS-0653121}

\address{127 Vincent Hall, University of Minnesota, Minneapolis, MN
55455} 
\email{krylov@math.umn.edu}


\begin{document}

\begin{abstract}
The solvability in $W^{2}_{p}(\bR^{d})$ spaces
is proved for second-order elliptic
equations with coefficients which are measurable
in one direction and VMO in the orthogonal directions
in each small ball with the direction depending on the ball.
This generalizes to a very large extent the case of equations
with
continuous or VMO coefficients.
\end{abstract}

\maketitle

\mysection{Introduction and main result}

In this article we are concerned with the
solvability in $W^{2}_{p}=
W^{2}_{p}(\bR^{d})$ of the 
equation
\begin{equation}
                       		                             \label{11.13.1}
Lu(x)-\lambda u(x) = f(x),
\end{equation}
where $L$ is a uniformly nondegenerate
elliptic differential operator
with bounded coefficients
 of the form
$$
Lu(x) = a^{ij}(x) u_{x^i x^j}(x) + b^{i}(x) u_{x^i}(x)
+ c(x) u(x)
$$ 
 in 
$$
   \bR^d
=\{x=(x^{1},...,x^{d}):x^{1},...,x^{d}\in\bR\}.
$$

We generalize the main result of
\cite{KK} where  the solvability
is established in the case
 that, roughly speaking, the coefficients $a^{ij}$ are measurable
with respect to $x^{1}$ and are in VMO with respect to 
$(x^{2},...,x^{d})$. Owing to a standard
localization procedure, this result admits an obvious
extension to the case in which for each ball $B\subset\bR^{d}$
of a {\em fixed\/} radius
 there exists a sufficiently regular diffeomorphism
that transforms equation \eqref{11.13.1} in $B$
into a similar equation with coefficients satisfying
the conditions of \cite{KK} in $B$. In particular,
one obtains the solvability if the matrix 
$a=(a^{ij})$ depends
on $|x|$ in a measurable way, is in VMO with
respect to the angular coordinates, and, say, is continuous at 
the origin.

The main goal of the present article is to 
show that in the above described generalization
 the radius of balls need not be fixed. 
{\color{black}In the end of this section we give an example
in which our result is applicable
contrarily to the result of \cite{KK}.}
 
We develop a new technique
which seems to be 
 applicable in many situations for elliptic and parabolic
equations with partially VMO coefficients as,
for instance, in \cite{DK}
and \cite{DKK}.
We only concentrate on elliptic equations in order to
make simpler the presentation of the method.
Generally, the theory of  elliptic 
 equations with partially VMO coefficients
is quite new and originated in \cite{KK} in contrast
with the case of completely VMO coefficients,
which appeared in \cite{CFL}, or the  
classical  
 case of
equations with continuous coefficients
treated in \cite{ADN}.
The reader can find
further references to articles and books
related to equations with VMO and partially VMO coefficients
in the above cited articles and the references therein.

In \cite{ADN} the main technical tool was the
theory of singular integrals, in particular, the
Calder\'on-Zygmund theorem. With development of Real
Analysis  later on
in many sources the theory of singular integrals
in applications to PDEs was replaced with using
the John-Nirenberg theorem or
Stampacchia interpolation theorem applied to
sharp functions. However, the theory of singular integrals
was used again in the paper \cite{CFL},
the results of which came as a real 
breakthrough in the theory of PDEs. Again later it turned out that
using singular integrals can be replaced with
appropriate other tools from Real Analysis such as
the Fefferman-Stein theorem. To the author it seems
highly unlikely that the theory of singular integrals
can be used to obtain 
even the main auxiliary result of \cite{KK}, which is
the basis of the present paper along with
  a new inequality of the Fefferman-Stein type
proved in Theorem \ref{theorem 7.1.1}.

In   connection with this new development
 it is instructive to recall
  that L. Bers and M. Schechter said in 1964
(see \cite{BJS}) that the linear theory of second order
elliptic PDEs ``is at present probably nearing
completion".

This paper deals with elliptic equations in nondivergence form.
A different technique is developed in several articles
by the authors of \cite{BL} for treating divergence type
equations. It would be interesting to know if
their methods could be applied to divergence
or nondivergence type equations
with coefficients satisfying our conditions. 
This could lead to extending our results to
equations in domains. So far we can only
deal with equations in the whole space or,
for that matter, with interior estimates.
Another restriction is that $p>2$.

Now we state our assumptions rigorously.

\begin{assumption}
                                           \label{assumption 6.30.1} 

The coefficients $a^{ij}$, $b^{i}$, and $c$
are measurable functions defined on $\bR^d$, $a^{ji} = a^{ij}$
for all $i,j=1,...,d$. There
exist positive constants $\delta\in(0,1)$ and $K$ such that
$$ 
|b^{i}(x)|  \le K,\quad i=1,...,d,
\qquad |c(x)| \leq K,
$$
$$
\delta |\xi|^2 \le   a^{ij}(x) \xi^i
\xi^j
\le\delta^{-1} |\xi|^2
$$ 
for any $x \in \bR^d$ and $\xi \in \bR^d$.
\end{assumption}

To state the second assumption denote by $A$ the set of
$d\times d$ symmetric matrix-valued
measurable functions $\bar{a}=
(\bar{a}^{ij}(t))$
of one variable
  $t\in\bR$ such that
$$
\delta |\xi|^2 \le   \bar{a}^{ij}(t) \xi^i
\xi^j
\le\delta^{-1} |\xi|^2
$$ 
for any $t\in \bR $ and $\xi \in \bR^d$.

Introduce $\Psi$ as 
the set of 
mappings $\psi:\bR^{d}\to \bR^{d}$ such that

(i) the mapping $\psi$ has an inverse $\psi^{-1}:
\bR^{d}\to \bR^{d}$;

(ii) the mappings $\psi$ and $\phi=\psi^{-1}$
are twice continuously differentiable and
$$
|\psi_{x}|+|\psi_{xx}|\leq\delta^{-1} ,
\quad
|\phi_{y}|+|\phi_{yy}|
\leq\delta^{-1}  .
$$

The following assumption contains a parameter $\gamma>0$,
which will be specified later. 
We denote  by $|B|$ the volume
of a Borel set $B\subset\bR^{d}$.

\begin{assumption}[$\gamma$]
                                           \label{assumption 6.30.2} 
There exists a  constant  $R_{0}>0$ such that
for any ball $B\subset\bR^{d}$ of radius less than $R_{0}$
one can find an
$\bar{a}\in A$ and a $\psi=(\psi^{1},...,\psi^{d})\in\Psi $ such that
\begin{equation}
                                                        \label{7.3.2}
\int_{B}|a(x)-\bar{a}(\psi^{1}(x))|\,dx\leq\gamma|B|.
\end{equation}
\end{assumption}

\begin{remark}
                                                  \label{remark 7.3.1}

Assumption \ref{assumption 6.30.2} ($\gamma$) is 
obviously satisfied with
any $\gamma>0$ if $a$ is uniformly continuous
as, for instance, in \cite{ADN}. If 
Assumption \ref{assumption 6.30.2}~($\gamma$) is 
  satisfied with
any $\gamma>0$ and {\em constant\/} $\bar{a}$
(perhaps, changing with $B$), then
one says that $a$ belongs to VMO. This case was
first treated in \cite{CFL}. In \cite{KK}
the solvability in $W^{2}_{p}$ was proved under
Assumption \ref{assumption 6.30.2}~($\gamma$)
with a fixed function 
{\color{black}$\psi$}, which
is not allowed to change with $B$. {\color{black}
(Actually, $\psi=x$ in \cite{KK}, but changing
coordinates shows that the result holds for any $\psi\in\Psi$.)}
{\color{black}By using partitions of unity}
the latter restriction on $\psi$ can
be  easily somewhat relaxed to allow mappings such that
in each ball $B$ 
 of  radius exactly $R_{0}$ there is a mapping $\psi$
which would suit all subballs inside $B$.
\end{remark}

As usual,  by $W_p^2=W_p^2(\bR^d)$ we mean the Sobolev
space  on $\bR^d$. Set   $\cL_p = \cL_p(\bR^d)$.

Here is our main result.

\begin{theorem}
                                              \label{theorem 6.30.1}
 Take a $p \in (2, \infty)$.
  Then there
exists a constant $\gamma=\gamma(d,\delta,p)>0$ such that
if Assumptions \ref{assumption 6.30.1} and
\ref{assumption 6.30.2} $(\gamma)$ are satisfied then  
for any $\lambda 
\geq\lambda_0(d,\delta,K,p,R_{0})\geq1$     and any
$f \in \cL_p$, there exists a unique $u \in W_p^2$ satisfying 
\eqref{11.13.1} in $\bR^{d}$.

Furthermore, there is a constant $N$, depending only on $d$, $\delta$,
$K$, $p$, and   $R_{0}$, such that, for any $\lambda \ge
\lambda_0$ and $u \in W_p^2$,
\begin{equation}         
                                                 \label{7.1.8}
\lambda \| u \|_{\cL_p} + \sqrt{\lambda} \| u_x \|_{\cL_p} + \| u_{xx}
\|_{\cL_p}
\le N \| Lu - \lambda u \|_{\cL_p}.
\end{equation}
\end{theorem}

The proof of this theorem is given in Section \ref{section 7.2.4}
after we prepare the necessary auxiliary results 
in Section \ref{section 7.2.3}, which in turn require
some general facts proved in Section \ref{section 7.2.2}.

We finish the section by 
{\color{black} giving the example we were talking about
above.
Let $f$ be a measurable function on
$\bR$ with support in
the interval $(1/2,1)$ and such that $|f|\leq1$.
Introduce $\xi(x)=\ln(|x|\wedge1)$, $x\in\bR$.
It is well known that $\xi\in BMO$. Then for $\varepsilon>0$
the function $\varepsilon\xi$ is also in $BMO$ and its $BMO$-norm
can be made as small as we like on the account of choosing
$\varepsilon$ small enough. The same is true
for $\eta=\sin(\varepsilon\xi)$ and $\zeta(x)=\eta(4x-3)$
with the latter function having support in the interval $(1/2,1)$.
Next take a large $\kappa\geq4$ and for real $x,y$ and $z=(x,y)$ introduce
$$
a(z)=\sum_{n=0}^{\infty}f(\kappa^{2n}x)\zeta(\kappa^{2n}y)
+\sum_{n=0}^{\infty}f(\kappa^{2n+1}y)\zeta(\kappa^{2n+1}x).
$$
Notice that  the support  of $f(\kappa^{r}\cdot)\zeta(\kappa^{r}\cdot)$
  belongs to $Q_{r}:=
(\kappa^{-r}/2,\kappa^{-r})^{2}$.

Now, for a square $Q=I
\times J\subset\bR^{2}$ we are going to
estimate the left-hand side of \eqref{7.3.2} with $Q$ and $z$ in place of
$B$ and $x$, respectively, and with $\psi$ equal to either $x$ or $y$.
For brevity we denote the modified left-hand side of \eqref{7.3.2}
by $M$.

Define $\tau$
as the least integer $k\geq0$ such that $Q\cap Q_{k}\ne\emptyset$.
If there are no such $k$'s, then $M=0$.
If $\tau$ is an even number we set $\psi=x$ and $\bar{a}(x)
=f(\kappa^{\tau}x)\bar{\zeta}$, where $\bar{\zeta}$
is the integral average of $\zeta(\kappa^{\tau}y)$ over $J$. Then 
\begin{equation}
                                                     \label{8.21.1}
M\leq\int_{I}|f(\kappa^{\tau}x)|\,dx\int_{J}|\zeta(\kappa^{\tau}
y)-\bar{\zeta}|\,dy
+\sum_{i=\tau+1}^{\infty}|Q\cap Q_{i}|.
\end{equation}

On the right, the first term is less than $|Q|\,|\zeta|_{BMO}$.
Also observe that if $i\geq\tau+1$
and $Q\cap Q_{i}\ne\emptyset$, then the lengths of $I$ and $J$
are at least $\kappa^{-\tau}/2-\kappa^{-i}$, which is larger
than $\kappa^{-\tau}/4$ since $\kappa\geq4$. Hence, in that case
$ 
|Q\cap Q_{i}|\leq|Q_{i}|=4^{-1}\kappa^{-2i}
\leq4\kappa^{2\tau-2i}|Q|
$ 
implying that the infinite sum in \eqref{8.21.1} is less than
$4(\kappa^{2}-1)^{-1}|Q|$. We see that in the case that $\tau$
is an even number   $M\leq\gamma|Q|$  with any fixed $\gamma>0$
provided that we choose sufficiently small $\varepsilon$
and sufficiently large $\kappa$.

In case $\tau$ is odd interchanging $x$ and $y$ leads to the same
conclusion and this easily shows
 that \eqref{7.3.2} holds indeed in its original form. Obviously,
functions like
the above $a$ cannot be treated by methods of \cite{KK}
even modified in the way outlined in 
Remark \ref{remark 7.3.1}.

}

The author wishes to thank Hongjie Dong for pointing out
several flaws in the first draft of the article.

\mysection{A partial version of the Fefferman-Stein theorem}
                                                 \label{section 7.2.2}

First we recall a few standard notions and facts
related to partitions and stopping times. All of them can be found in
many books; we follow the exposition in  \cite{Kr08}.

Let $(\Omega,\cF,\mu)$ be a complete measure space with
a $\sigma$-finite measure $\mu$, such that
$$
\mu(\Omega)=\infty.
$$
Let $\cF^{0}$ be the subset
 of $\cF$ consisting of all                                    
sets $A$ such that $\mu(A)<\infty$.
For $p\in[1,\infty)$ set  $\cL_{p}(\Omega)=\cL_{p}(\Omega,\cF,\mu)$.
By $\cL^{0}$ we denote a fixed dense subset
of $\cL_{1}(\Omega)$.
For any $A\in\cF$ we 
set
$$
|A|=\mu(A).
$$
For $A\in\cF^{0}$ and   functions $f$ 
summable on $A$ we 
  use the notation 
$$
f_{A}=\dashint_{A}f\,\mu(dx):=\frac{1}{|A|}\int_{A}f(x)\,\mu(dx)  
\quad \bigg(\frac{0}{0}:=0\bigg)
$$
 for  the average value of $f$ over $A$.

\begin{definition}
                                          \label{definition 5.30.1}
Let $\bZ=\{n: n=0,\pm1,\pm2,...\}$ and let
 $(\bC_{n},n\in\bZ)$ be
a sequence of   partitions of $\Omega$ each consisting of  
countably many disjoint
  sets $C\in\bC_{n} $ and such that $\bC_{n}\subset
\cF^{0}$ for each $n$. 
For each $x\in\Omega$ and $n\in \bZ$ there exists (a unique)
$C\in\bC_{n}$ such that $x\in C$. We denote 
this $C$ by
$C_{n}(x)$.

The sequence $(\bC_{n},n\in\bZ)$ is called   a filtration
of partitions 
if the following conditions are satisfied.

(i) The elements of partitions 
are ``large"     for big negative $n$'s and ``small"  
for  big positive
$n$'s:
$$
\inf_{C\in\bC_{n}}|C|\to\infty\quad\hbox{as}\quad n\to-\infty,\quad
 \lim_{n\to\infty}f_{C_{n}(x)}=f(x)\quad\hbox{(a.e.)}\quad
\forall f\in\cL^{0}.
 $$

(ii) The partitions are nested: for each $n$ and $C\in\bC_{n}$
there is a (unique) $C'\in\bC_{n-1}$
 such that $C\subset C'$.

(iii) The following regularity property holds: for any $n$, $C$, and $C'$
as in (ii) we have
$$
|C'|\leq N_{0}|C|,
$$
 where $N_{0}$
is a constant independent of $n,C, C'$.
\end{definition}

Observe that since the elements of partition $\bC_{n}$
become large as $n\to-\infty$, we have $N_{0}>1$.

The only example of a filtration of partitions 
important for this article
in the case that $\Omega=\bR^{d }$ with Lebesgue measure $\mu$
is given
by   dyadic cubes, that is, by 
$$
\bC_{n}=\{C_{n}(i_{1},...,i_{d}),
i_{1},...,i_{d}\in\bZ\},
$$
 where
$$
C_{n}(i_{1},...,i_{d})=[i_{1}2^{-n},(i_{1}+1)2^{-n})\times...\times
[i_{d}2^{-n},(i_{d}+1)2^{-n}).
$$
In this case, to satisfy requirement (i) in Definition
\ref{definition 5.30.1}, one can take $\cL^{0}$ as the set of
continuous functions with compact support.

\begin{definition}
                                          \label{definition 5.30.2}
 Let $\bC_{n}$, $n\in\bZ$, be a filtration of partitions
of $\Omega$.

(i) Let $\tau=\tau(x)$ be a function on $\Omega$ with values in
$\{\infty,0,\pm1,\pm2,...\}$. The function $\tau$
is called  a stopping time 
(relative to the filtration)
if, for each $n=0,\pm1,\pm2,...$, the set 
$$
\{x:\tau(x)=n\}
$$ is either empty or else is the union of
some elements of $\bC_{n}$.

(ii) For a function $f\in \cL_{1 }(\Omega)$ and
$n\in\bZ$,
we denote
$$
f_{|n}(x)= \dashint_{C_{n}(x)}f(y)\,\mu(dy).
$$
We read $f_{|n}$ as ``$f$ given $\bC_{n}$'', continuing to
borrow  the terminology from probability theory.
If we are also given a stopping time $\tau$, we 
let
$$
f_{|\tau}(x) =f_{|\tau(x)}(x)
$$
 for those $x$ for which
$\tau(x)<\infty$  and $f_{|\tau}(x)=f(x)$ otherwise.
\end{definition}

The simplest example of a stopping time is given by
$\tau(x)\equiv0$. It is also known that
if $g\in\cL_{1}(\Omega)$ and a constant $\lambda>0$, then
$$
\tau(x)=\inf\{n\in\bZ:g_{|n}(x)>\lambda\}\quad(\inf\emptyset:=
\infty)
$$
is a stopping time and if, in addition, $g\geq0$, then
$g_{|\tau}\leq N_{0}\lambda$ (a.e.).

For $f\in\cL_{1}(\Omega)$ we denote
$$
\cM f=\sup_{n\in\bZ}|f|_{|n }.
$$
It is known that for any $f\in\cL_{1}(\Omega)$ and $p\in(1,\infty)$
\begin{equation}
                                                \label{7.1.6}
\|\cM f\|_{\cL_{p}(\Omega)}\leq q\|f\|_{\cL_{p}(\Omega)},
\end{equation}
where $q=p/(p-1)$.

In the remaining part of the section we consider two functions
$u,v\in\cL_{1}(\Omega)$ and a nonnegative measurable function $g$
on $\Omega$. {\color{black}
\begin{lemma}
                                                 \label{lemma 10.27.1}
Assume that $0\leq u\leq v$  and for any $n\in\bZ$
and
$C\in\bC_{n}$ we have
\begin{equation}
                                                      \label{10.27.2}
\int_{C}(u  -v _{C})_{+}\,\mu(dx)\leq \int_{C}g(x)\,\mu(dx).
\end{equation}
Then for any $\lambda>0$
\begin{equation}
                                                      \label{10.27.1}
 |\{x: u(x) \geq\lambda\}|\leq 2\lambda^{-1}
\int_{\Omega}g(x)I_{\cM v(x)>\alpha \lambda}\,\mu(dx),
\end{equation}
where $\alpha=(2N_{0})^{-1}$.
\end{lemma}

Proof. Fix a $\lambda>0$ and define 
$$
\tau(x)=\inf\{n\in\bZ:v_{|n}(x)>\alpha\lambda\}.
$$
We know that $\tau$ is a stopping time and if $\tau(x)<\infty$,
then 
$$
  v_{|n}(x)\leq \lambda/2,\quad\forall n\leq\tau(x).
$$
We also know that $v_{|n}\to v$ (a.e.) as $n\to\infty$.
It follows that (a.e.)
$$
\{x:u(x)\geq\lambda\}=\{x:u(x)\geq\lambda,\tau(x)<\infty\}
$$
$$
=\{x:u(x)\geq\lambda,  v_{|\tau}\leq \lambda/2\}
=\bigcup_{n\in\bZ}\bigcup_{C\in \bC^{\tau}_{n}}A_{n}(C),
$$
where
$$
A_{n}(C):=\{x\in C:u(x)\geq\lambda, v_{|n}\leq \lambda/2\},
$$
and $\bC^{\tau}_{n}$ is the family of disjoint elements
of $\bC_{n}$ such that
$$
\{x:\tau(x)=n\}=\bigcup_{C\in \bC^{\tau}_{n}}C.
$$

Next, for each $n\in\bZ$ and $C\in\bC_{n}$ on the set $A_{n}(C)$,
if it is not empty,
we have $v_{|n}=v_{C}$ and $u-v_{C}\geq\lambda/2$, so that
by Chebyshev's inequality and assumption \eqref{10.27.2}
$$
|A_{n}(C)|\leq 
2\lambda^{-1}\int_{C}g \,\mu(dx),
$$
$$
|\{x:u(x)\geq\lambda\}|\leq2\lambda^{-1}
\sum_{n\in\bZ}\sum_{C\in \bC^{\tau}_{n}}\int_{C}g \,\mu(dx)
=2\lambda^{-1}\int_{\Omega}gI_{\tau<\infty}\,\mu(dx).
$$
It only remains to observe that $\{\tau<\infty\}=\{\cM
v>\alpha\lambda\}$. The lemma is proved.
\begin{remark}
Obviously, the conditions of Lemma  \ref{lemma 10.27.1}
are satisfied with $g=(1/2)v^{\sharp}$ if $u=v$.
One of nice features of the lemma
is that under its conditions, for any measurable function
$a$ such that $1\leq a \leq2$, the functions $au$, $2v$,
and $2g$ also satisfy its conditions.
\end{remark}

To give conditions to verify assumption \eqref{10.27.2}
which are convenient in this
article,} we need the following.

\begin{assumption}
                                         \label{assumption 7.1.1}
We have $|u|\leq v$
and 
for any $n\in\bZ$ and $C\in\bC_{n}$ there exists a
measurable function  $u^{C}$ given on $C$ such that
$|u|\leq u^{C}\leq v$ on $C$ and
\begin{equation}
                                                      \label{6.29.3}
\big(\int_{C}|u -u _{C}|\,\mu(dx)\big)
\wedge\big(\int_{C}|u^{C}-u^{C}_{C}|\,\mu(dx)
\big)\leq\int_{C}g(x)\,\mu(dx).
\end{equation}
\end{assumption}

\begin{lemma}
Under Assumption \ref{assumption 7.1.1}   
for any $\lambda>0$ we have
\begin{equation}
                                                      \label{6.29.1}
 |\{x:|u(x)|\geq\lambda\}|\leq {\color{black}2}\lambda^{-1}
\int_{\Omega}g(x)I_{\cM v(x)>\alpha \lambda}\,\mu(dx),
\end{equation}
where $\alpha=(2N_{0})^{-1}$. Moreover if $u\geq0$, then
one can replace ${\color{black}2}\lambda^{-1}$ in \eqref{6.29.1} with
$ {\color{black} \lambda^{-1}}$.
 
\end{lemma}

Proof. First assume that $u \geq0$. 
{\color{black} Take an $n\in\bZ$ and a $C\in\bC_{n}$. If
$$
\int_{C}|u -u _{C}|\,\mu(dx)\leq\int_{C}g(x)\,\mu(dx).
$$
then, since $u\leq v$, we have $u_{C}\leq v_{C}$ and 
$$
(u -u _{C})+|u -u _{C}|=2(u -u _{C})_{+} \geq2
(u -v_{C})_{+},
$$
implying that \eqref{10.27.2} is satisfied with $g/2$ in place
of $g$. In case that
$$
\int_{C}|u^{C}-u^{C}_{C}|\,\mu(dx)\leq\int_{C}g(x)\,\mu(dx)
$$
we observe that $u^{C}\geq u $, $u^{C}_{C}\leq v_{C}$, so that
$$
(u^{C}-u^{C}_{C})+|u^{C}-u^{C}_{C}|=2(u^{C} -u^{C} _{C})_{+} \geq2
(u -v_{C})_{+},
$$
which again implies that \eqref{10.27.2} is satisfied with $g/2$ in
place of $g$.
}

In the general case we need only show that condition \eqref{6.29.3}
is almost preserved if we take $|u |$ in place of $u $. However,
for any   measurable set $C$ we have
$$
 \dashint_{C}\big|\,|u(x)|- |u|_{C} \big|\,\mu(dx)
=\dashint_{C}\big|\dashint_{C}(|u(x)|-|u(y)|)\,\mu(dy)\big|\,\mu(dx)
$$
\begin{equation}
                                                        \label{7.1.4}
\leq\dashint_{C} \dashint_{C} |u(x) - u(y)| \,\mu(dy) \,\mu(dx)
\leq 2\dashint_{C} |u(x) -c| \,\mu(dx),
\end{equation}
where $c$ is any constant. If we take $c=u_{C}$, then we see that
$|u |$ satisfies \eqref{6.29.3} with $2g$ in place of $g$.
The lemma is proved.

Now we are ready to prove a partial version of the Fefferman-Stein
theorem about sharp functions.

\begin{theorem}
                                          \label{theorem 7.1.1}
Under Assumption \ref{assumption 7.1.1}   
for any $p\in(1,\infty)$ we have
\begin{equation}
                                                      \label{7.1.5}
\|u\|_{L_p(\Omega)}^{p}\leq N(p,N_{0})\|g\|_{L_p(\Omega)} 
\|v\|_{L_p(\Omega)}^{p-1}.
\end{equation}
{\color{black} The same conclusion holds
under the assumptions of Lemma \ref{lemma 10.27.1}.}
 
\end{theorem}

Proof. We have
$$
\|u\|_{L_p(\Omega)}^{p}=
\int_{0}^{\infty} |\{x:|u(x)|\geq\lambda^{1/p}\}|\,d\lambda
$$
$$
\leq {\color{black}2}
\int_{\Omega}g(x)\big(
\int_{0}^{\infty}\lambda^{-1/p}I_{\cM v(x)>\alpha \lambda^{1/p}}
\,d\lambda\big)\,\mu(dx)
$$
$$
={\color{black}2}q\alpha^{1-p}\int_{\Omega}g (\cM v)^{p-1}\,\mu(dx),
$$
where $q=p/(p-1)$. By using H\"older's inequality
and \eqref{7.1.6}, we come to
\eqref{7.1.5}. The theorem is proved.

\begin{remark}
In the dyadic version of the
original Fefferman-Stein theorem $u^{C}=u$, $v=|u|$,
and $g$ is the sharp function $u^{\sharp}$ of $u$. In that case,
assuming that $u\in \cL_{p}(\Omega)$, we get from
\eqref{7.1.5} the Fefferman-Stein inequality
$\|u\|_{L_p(\Omega)} \leq N \|u^{\sharp}\|_{L_p(\Omega)}$.
\end{remark}

\mysection{Auxiliary results}
                                                 \label{section 7.2.3}

We denote by $B_{r}(x)$ the open ball in $\bR^{d}$
of radius $r$
centered at $x$. Set $B_{r}=B_{r}(0)$ and introduce
$\bB$ as the family of balls in $\bR^{d}$.
For a Borel set $B\subset\bR^{d}$ of nonzero Lebesgue measure
and a measurable function $f$ we 
define
$$
f_{B}:=\dashint_{B}f(x)\,dx:=\frac{1}{|B|}\int_{B}f(x)\,dx,
$$
whenever the last integral is finite.
The following is  Lemma 4.8
of \cite{KK}.

\begin{lemma}
                                          \label{lemma 6.30.1}  
Take an $\bar{a}\in A$ and set
\begin{equation}         
                                                 \label{6.30.6}
\bar{L} u(x)=\bar{a}^{ij}( x^{1})u_{x^{i}x^{j}}(x).
\end{equation}
There exists a
constant  $N = N(d,\delta)$  such that, for any $\kappa \ge 4$, $r
> 0$,   $u \in C_0^{\infty}(\bR^d)$, and $i,j\in\{1,...,d\}$
satisfying $ij>1$ we have
$$
\dashint_{B_r} |u_{ x^{i}x^{j}} - 
\left(u_{x^{i}x^{j}}\right)_{B_r}|^2 \, dx 
\le  N \kappa^{d } \left( |\bar{L} u|^2 \right)_{B_{\kappa r}}   + N
\kappa^{-2}  
\left( |u_{xx}|^{2  } \right)_{B_{\kappa r}}  .
$$

\end{lemma}

We need a version of this lemma 
for operators of a more general form.

\begin{lemma}
                                          \label{lemma 7.3.2}  
Take an $\bar{a}\in A$ and a $\psi\in\Psi$ and set
\begin{equation}         
                                                 \label{7.3.5}
\bar{L} u(x)=\bar{a}^{kn}( y^{1})
\phi^{i}_{y^{n}}(y)\phi^{j}_{y^{m}}(y)u_{x^{i}x^{j}}(x),
\end{equation}
where $y=\psi(x)$ and $\phi=\psi^{-1}$. Then
there exist   
constants  $N = N(d,\delta)$ 
and $\chi = \chi(d,\delta)\geq1$ such that, for any $\kappa \ge 4$, $r
> 0$,   $u \in C_0^{\infty}(\bR^d)$, and $i,j\in\{1,...,d\}$
satisfying $ij>1$ we have
$$
\dashint_{B_r} |u_{ij} - 
\left(u_{ij}\right)_{B_r}|^2 \, dx 
\le  N \kappa^{d } \left( |\bar{L} u|^2 \right)_{B_{\chi\kappa r}}
$$
\begin{equation}         
                                                 \label{6.30.5}  
+N \kappa^{d } \left( |u_{x}|^2 \right)_{B_{\chi\kappa r}} 
 + N
\kappa^{-2}  
\left( |u_{xx}|^{2  } \right)_{B_{\chi\kappa r}}  ,
\end{equation}
where $u_{ij}(x)$ are defined by
\begin{equation}
                                                          \label{7.4.1}
u_{ij}( \phi(y))=v_{y^{i}y^{j}}(y),\quad v(y)=u( \phi(y)),
\quad \phi=\psi^{-1}.
\end{equation}

\end{lemma}

Proof. Without loss of generality we assume that $\psi(0)=0$.
  Also set $f=\bar{L} u$
and observe that
\begin{equation}
                                                          \label{7.3.6}
\bar{a}^{kn}(y^{1})v_{y^{k}y^{n}}(y)+\tilde{b}^{k}(y)v_{y^{k}}(y)
=f(\phi(y)),
\end{equation}
 where  
$$
 \tilde{b}^{k}(y)=
\bar{a}^{kn}( y^{1})
\phi^{i}_{y^{n}}(y)\phi^{j}_{y^{m}}(y)\psi^{k}_{x^{i}x^{j}}(x),\quad
x=\phi(y).
$$

Next we apply Lemma \ref{lemma 6.30.1}  to the operator
$$
\bar{L}_{y}v(y)=\bar{a}^{kn}(y^{1})v_{y^{k}y^{n}}(y)
$$
and for any $\rho>0$ find
\begin{equation}         
                                                 \label{7.3.7}
\dashint_{B_\rho} |v_{ y^{i}y^{j}} - 
\left(v_{y^{i}y^{j}}\right)_{B_\rho}|^2 \, dy 
\le  N \kappa^{d } \left( |\bar{L}_{y} v|^2 \right)_{B_{\kappa \rho}}   + N
\kappa^{-2}  
\left( |v_{yy}|^{2  } \right)_{B_{\kappa \rho}}  .
\end{equation}

To transform this inequality we use the simple observation
that there exist   constants  $N,\chi<\infty$ depending only on
$d$ and $\delta$ such that for any nonnegative measurable
function $g$ we have
$$
\dashint_{B_{\rho}}f(x)\,dx\leq N
\dashint_{B_{\rho\sqrt{\chi}}}f(\phi(y))\,dy,
\quad
\dashint_{B_{\rho}}f(\phi(y))\,dy\leq
N\dashint_{B_{ \rho\sqrt{\chi}}}f(x)\,dx.
$$
Using this and closely following \eqref{7.1.4} we find
$$
\dashint_{B_r} |u_{ij} - 
\left(u_{ij}\right)_{B_r}|^2 \, dx\leq
\dashint_{B_r}\dashint_{B_r}|u_{ij}(x_{1})-u_{ij}(x_{2}|^{2}
\,dx_{1}dx_{2}
$$
$$
\leq N
\dashint_{B_{  r\sqrt{\chi}}}\dashint_{B_{  r\sqrt{\chi}}}
|v_{y^{i}y^{j}}(y_{1})-v_{y^{i}y^{j}}(y_{2})|^{2}
\,dy_{1}dy_{2}
$$
$$
\leq N
\dashint_{B_{  r\sqrt{\chi}}} |v_{ y^{i}y^{j}} - 
\left(v_{y^{i}y^{j}}\right)_{B_\rho}|^2 \, dy .
$$

Furthermore, for $y=\psi(x)$ 
obviously $|v_{yy}(y)|\leq N(|u_{xx}(x)|+|u_{x}(x)|)$
and by \eqref{7.3.6} also $|\bar{L}_{y} v(y)|\leq|\bar{L}u(x)|
+N|u_{x}(x)|$. By combining the above observations we 
immediately obtain \eqref{6.30.5} from \eqref{7.3.7}.
The lemma is proved.
 
Set
$$
L_{0}u(x)=a^{ij}( x)u_{x^{i}x^{j}}(x).
$$
In the following lemma we prepare to check Assumption 
\ref{assumption 7.1.1} for some functions to be introduced later
and closely related to $u_{ij}$.
However, we still have $B_{r}$ in place of~$C$.

\begin{lemma}
                                          \label{lemma 6.30.2}  
(i) Suppose that Assumptions   \ref{assumption 6.30.1} and
\ref{assumption 6.30.2} $(\gamma)$ are satisfied.

(ii) Let  $\mu$, $\nu \in (1,\infty)$, $\kappa \ge 4$,  
  and $r>0$
be some numbers such that $1/\mu + 1/\nu =
1$. 

Then
there exist  a mapping $\psi\in\Psi $ and  
constants  $N = N(d,\delta,\mu)$ and $\chi = \chi(d,\delta )\geq1$
 such that, for any   
  $ C_0^{\infty} $ function $u$, vanishing outside
  a ball of radius $R\leq R_{0}$, and  $i,j\in\{1,...,d\}$
satisfying $ij>1$ we have
$$
\dashint_{ B_r } |u_{ij} - 
\left(u_{ij} \right)_{B_r}|^2 \, dx 
\le  N \kappa^{d } \left( |L_{0}u|^2 \right)_{B_{\chi \kappa r}} 
+ N \kappa^{d } \left( |u_{x}|^2 \right)_{B_{\chi \kappa r}} 
$$
\begin{equation}         
                                                 \label{6.30.8}
  + N( 
\kappa^{d } 
R^{2}+
\kappa^{-2})  
\left( |u_{xx}|^{2  } \right)_{B_{\chi \kappa r}} + 
 N\kappa^{d }
\gamma^{1/\nu}\left( |u_{xx}|^{2\mu} 
\right)^{1/\mu}_{B_{\chi\kappa r}},
\end{equation}
where $u_{ij}(x)$ are defined by \eqref{7.4.1}.

\end{lemma}

Proof.  We take $\chi$ from Lemma
\ref{lemma 7.3.2}   and split  the proof into
two parts.

{\em Case $\chi \kappa r<R$\/}. Take
a $\psi\in\Psi$ and an $\hat{a}\in A$ such that
\begin{equation}         
                                                 \label{6.30.7}
\dashint_{B_{\chi\kappa r}}|a(x)-\hat{a}(\psi^{1}(x))|\,dx\leq\gamma.
\end{equation}
Reducing $\delta$ if necessary we may assume that,
for an $\bar{a}\in\ A$, we have
\begin{equation}         
                                                 \label{7.4.5} 
\hat{a}^{ij}(t)=\bar{a}^{kn}(t)
\phi^{i}_{y^{n}}(y_{0})\phi^{j}_{y^{m}}(y_{0}).
\end{equation}
where $y_{0}=\psi(0)$.
Then introduce $\bar{L}$ by \eqref{7.3.5} and set
$$
\hat{L}u(x)=\hat{a}^{ij}(\psi^{1}(x))u_{x^{i}x^{j}}(x).
$$

Observe that for $y=\psi(x)$ and $|x|\leq\chi \kappa r$ we have
$|y-y_{0}|\leq N(d,\delta)\chi \kappa r$ and
$$
|(\bar{L}-\hat{L})u(x)|=\big|
\bar{a}^{kn}( y^{1})\big(
\phi^{i}_{y^{n}}(y)\phi^{j}_{y^{m}}(y)- 
\phi^{i}_{y^{n}}(y_{0})\phi^{j}_{y^{m}}(y_{0})\big)u_{x^{i}x^{j}}(x)\big|
$$
\begin{equation}
                                                          \label{7.4.2}
  \leq NR|u_{xx}(x)|.
\end{equation}
This and \eqref{6.30.5} yield
$$
\dashint_{B_r} |u_{ij} - 
\left(u_{ij}\right)_{B_r}|^2 \, dx 
\le  N \kappa^{d } \left( |\hat{L} u|^2 \right)_{B_{\chi\kappa r}}
+ N \kappa^{d } R^{2}\left( |u_{xx}|^2 \right)_{B_{\chi\kappa r}}
$$
\begin{equation}         
                                                 \label{7.4.3}  
+N \kappa^{d } \left( |u_{x}|^2 \right)_{B_{\chi\kappa r}} 
 + N
\kappa^{-2}  
\left( |u_{xx}|^{2  } \right)_{B_{\chi\kappa r}}  .
\end{equation}
After that it only remains to notice that
$$
\left( |\hat{L} u|^2 \right)_{B_{\chi\kappa r}}
\leq 2\left( |L_{0}u|^2 \right)_{B_{\chi\kappa r}}
+2\left( |(\hat{L} -L_{0}) u|^2 \right)_{B_{\chi\kappa r}}.
$$
and by  H\"older's inequality and \eqref{6.30.7}
\begin{equation}         
                                                 \label{6.30.9}
\left( |(\hat{L}-L_{0}) u|^2 \right)_{B_{\chi\kappa r}}
\leq N\left( |u_{xx}|^{2\mu} \right)^{1/\mu}_{B_{\chi\kappa r}}
\gamma^{1/\nu},
\end{equation}
which yields \eqref{6.30.8}.

{\em Case $\chi\kappa r\geq R$\/}. Let $u=0$ outside $B_{R}(x_{0})$.
Take
a $\psi\in\Psi$ and an $\hat{a}\in A$ such that
$$
\dashint_{B_{R}(x_{0})}|a(x)-\hat{a}(\psi^{1}(x))|\,dx\leq\gamma,
$$
 define $\bar{a} $ by \eqref{7.4.5} 
with $y_{0}=\psi(x_{0})$, and define $\hat{L}$ and $\bar{L}$
as above. Then on the support of $u$ we still have \eqref{7.4.2}
and hence \eqref{7.4.3}  holds again.
Finally,
$$
\left( |(\hat{L}-L_{0}) u|^2 \right)_{B_{\chi\kappa r}}
=\left(I_{B_{R}(x_{0})} |(\hat{L}-L_{0}) u|^2 \right)_{B_{\chi\kappa r}}
$$
$$
\leq N\left( |u_{xx}|^{2\mu} \right)^{1/\mu}_{B_{\chi\kappa r}}J,
$$
where
$$
J^{\nu} :=\frac{1}{|B_{\chi\kappa r}|}
\int_{B_{\chi\kappa r}\cap B_{R}(x_{0})}
|a(x)-\hat{a}(\psi^{1}(x))|\,dx
$$
$$
\leq\frac{1}{| B_{R}(x_{0})|}
\int_{ B_{R}(x_{0})}
|a(x)-\hat{a}(\psi^{1}(x))|\,dx\leq\gamma.
$$
It is seen that \eqref{6.30.9} is true again and the lemma is proved.

In the next lemma by $\bC_{n},n\in\bZ$, we mean the filtration of
dyadic cubes in $\bR^{d}$ and by $\bM f$ the classical
maximal function of $f$ defined by
$$
\bM f(x)=\sup_{B\in\bB:B\ni x}\dashint_{B}|f(y)|\,dy.
$$

\begin{lemma}
                                          \label{lemma 7.1.1} 
(i) Suppose that Assumptions   \ref{assumption 6.30.1} and
\ref{assumption 6.30.2} $(\gamma)$ are satisfied.

(ii) Let  $\mu$, $\nu \in (1,\infty)$, and $\kappa \ge 4$ 
be some numbers such that $1/\mu + 1/\nu =
1$. 

Then
for any $n\in\bZ$ and $C\in\bC_{n}$
there exist    a mapping $\psi\in\Psi $  and a
constant  $N = N(d,\delta,\mu)$  such that, for any   
  $ C_0^{\infty} $ function $u$, vanishing outside
  a ball of radius $R\leq R_{0}$,  and  $i,j\in\{1,...,d\}$
satisfying $ij>1$ we have
\begin{equation}         
                                                 \label{7.1.1}
\int_{C} |u_{ij} - 
 (u_{ij} )_{C}|  \, dx 
\leq N\int_{C}g\,dx,
\end{equation}
where  $u_{ij}(x)$ are defined by \eqref{7.4.1} and
 $g$ is a nonnegative function satisfying
$$
g^{2}=  \kappa^{d }( \bM( |L_{0}u|^2 ) +\bM( |u_{x}|^2 ))
$$
$$
 + ( 
\kappa^{d }  
R^{2}+\kappa^{-2})  
\bM( |u_{xx}|^{2  } )  + 
 \kappa^{d }
\gamma^{1/\nu}\left(\bM( |u_{xx}|^{2\mu} 
) 
\right)^{1/\mu}.
$$
Furthermore,  
\begin{equation}         
                                                 \label{7.3.4}
|u_{xx}|\leq N\sum_{ij>1}|u_{ij}|+N|u_{x}|+
N|L_{0}u|.
\end{equation}

\end{lemma}

Proof. Let $B$ be the smallest ball containing $C$
and let $B'$ be the concentric ball of radius $\chi\kappa r$, where
$r$ is the radius of $B$ and $\chi$ is taken from
Lemma \ref{lemma 6.30.2} . One can certainly shift the origin in
the situation of Lemma \ref{lemma 6.30.2} and hence 
for $ij>1$ and an appropriate $\psi\in\Psi$   
$$
\dashint_{ B } |u_{ij} - 
\left(u_{ij} \right)_{B }|^2 \, dx 
\le  N_{1} \kappa^{d } \left( |L_{0}u|^2 \right)_{B'} 
+ N_{1} \kappa^{d } \left( |u_{x}|^2 \right)_{B'} 
$$
\begin{equation}         
                                                 \label{7.4.05}
  + N_{1}({\color{black}\kappa^{d}}R^{2}+
\kappa^{-2})  
\left( |u_{xx}|^{2  } \right)_{B'} + 
 N_{1}\kappa^{d }
\gamma^{1/\nu}\left( |u_{xx}|^{2\mu} 
\right)^{1/\mu}_{B'},
\end{equation}
where $N_{1}=N(d,\delta,\mu)$. Obviously, the right-hand side of
\eqref{7.4.05} is less than $N_{1}g^{2}(x)$ for any $x\in C$
(and for that matter, for any $x\in B'$).
In particular, the square root of the right-hand side of
\eqref{7.4.05} is less than
$$
N_{1}^{1/2}\dashint_{C}g\,dx.
$$
After that, to finish proving the first
assertion of the lemma, it only remains to use H\"older's inequality
showing that
$$
J:=\dashint_{B } |u_{ij} - 
\left(u_{ij}\right)_{B }| \, dx
\leq\big(\dashint_{B } |u_{ij} - 
\left(u_{ij}\right)_{B }|^2 \, dx\big)^{1/2} 
$$
and observe that
$$
\dashint_{C} |u_{ij} - 
\left(u_{ij}\right)_{C }|  \, dx\leq
\dashint_{C}\dashint_{C} |u_{ij}(x) - 
 u_{ij}(y)|  \, dxdy
$$
$$
\leq N(d)
\dashint_{B}\dashint_{B} |u_{ij}(x) - 
 u_{ij}(y)|  \, dxdy\leq NJ.
$$

To prove the second assertion,
define $f=L_{0}u$, $v(\psi(x))=u(x)$, and by changing variables
introduce an operator $\hat{L}$ such that $\hat{L}v(y)=f(\phi(y))$.
Then
$$
|v_{yy}|\leq N\sum_{ij>1}|v_{y^{i}y^{j}}|
+N|\hat{L}v|+N|v_{y}|.
$$
By adding to this that $|u_{xx}(x)|\leq N|v_{yy}(y)|+N|u_{x}(x)|$
for $y=\psi(x)$,
we come to \eqref{7.3.4}.
The lemma is proved.

\begin{lemma}
                                          \label{lemma 7.1.2} 
Let $p\in(2,\infty)$. We assert that there exist    constants
$\gamma=\gamma(d,\delta,p)>0$ and $R=R(d,\delta,p,R_{0})\in(0,R_{0}]$
 such that if 
Assumptions   \ref{assumption 6.30.1} and
\ref{assumption 6.30.2} $(\gamma)$ are satisfied,
then for any $C^{\infty}_{0}$ function $u$ vanishing
outside   a ball of radius $R$ we have
\begin{equation} 
\|u_{xx}\|_{\cL_{p}}\leq N(\|L_{0}u\|_{\cL_{p}}  \label{7.1.3}
+\|u_{x}\|_{\cL_{p}}),
\end{equation}
where $N=N(d,\delta,p)$.

\end{lemma}

Proof. 
For the moment we suppose that Assumptions   \ref{assumption 6.30.1} and
\ref{assumption 6.30.2}~$(\gamma)$ are satisfied with a 
constant $\gamma>0$
and will choose it appropriately near the end of the proof.

Take a number $\kappa\geq4$ and set $\mu=(2+p)/4$ ($\mu>1,2\mu<p$).
Also take an $n\in\bZ$ and a $C\in\bC_{n}$ and take a $\psi\in\Psi$
from Lemma \ref{lemma 7.1.1}.
Finally, take a $C^{\infty}_{0}$ function $u$ vanishing
outside   a ball of radius $R$, introduce $u_{ij}$ by
\eqref{7.4.1}, and set 
$$
L_{0}u=f,\quad 
U=|u_{xx}|,\quad  
U  
^{C}=\sum_{ij>1}|u_{ij}|+|u_{x}|+|f|,
\quad
V=|u_{xx}|+|u_{x}|+|f|.
$$ 

We want to apply Theorem \ref{theorem 7.1.1}.
Estimate \eqref{7.3.4} says that
 $
U\leq NU^{C} 
 $.
Furthermore, obviously $U^{C}\leq NV$.
Also, similarly to 
\eqref{7.1.4}
$$
\dashint_{C}|U^{C}-U^{C}_{C}|\,dx\leq 2
\sum_{ij>1}\dashint_{C}|u_{ij}-(u_{ij})_{C}|\,dx
$$
$$
+2\dashint_{C}|u_{x}-(u_{x})_{C}|\,dx+2\dashint_{C}|f-f_{C}|\,dx.
$$
We estimate the sum over $ij>1$ by using Lemma \ref{lemma 7.1.1} 
and observe that
$$
\dashint_{C}|f-f_{C}|\,dx\leq 2|f|_{C}\leq2\bM f(x)\quad\forall
x\in C,
$$
$$
 \dashint_{C}|f-f_{C}|\,dx
\leq2\dashint_{C}\bM f\,dx,\quad
\dashint_{C}|u_{x}-(u_{x})_{C}|\,dx\leq 2\dashint_{C}\bM|u_{x}|\,dx.
$$
Hence
$$
\dashint_{C}|U^{C}-U^{C}_{C}|\,dx\leq N
\dashint_{C}(g+\bM|u_{x}|+\bM f)\,dx,
$$
where $g$ is defined in Lemma \ref{lemma 7.1.1} .

Since this holds for any $n\in\bZ$ and any  $C\in\bC_{n}$,
by Theorem \ref{theorem 7.1.1} we conclude
$$
\|u_{xx}\|_{\cL_{p}}=\|U\|_{\cL_{p}}
\leq N\|g+\bM|u_{x}|+\bM f\|_{\cL_{p}}^{1/p}\|V\|^{(p-1)/p}_{\cL_{p}}.
$$
By observing that
$$
\|V\| _{\cL_{p}}\leq\|u_{xx}\| _{\cL_{p}}+\|u_{x}\| _{\cL_{p}}+\|f\| _{\cL_{p}}
$$
and by Young's inequality 
$$
a^{1/p}b^{(p-1)/p}\leq
N(\varepsilon,p)a+\varepsilon b,\quad\forall
a,b,\varepsilon>0,
$$
 we 
easily get that
$$
\|u_{xx}\|_{\cL_{p}}
\leq N\|g+\bM|u_{x}|+\bM f\|_{\cL_{p}}+\|u_{x}\| _{\cL_{p}}+\|f\| _{\cL_{p}} .
$$
Next, by applying the Hardy-Littlewood maximal function theorem
and using the fact that $p/(2\mu)>1$ and $p>2$ we find
$$
\|u_{xx}\|_{\cL_{p}}\leq N_{1} \kappa^{d/2}\|f\|_{\cL_{p}}
+N_{1} \kappa^{d/2}\|u_{x}\|_{\cL_{p}}
$$
$$
+N_{1}( 
\kappa^{d/2}  
R+\kappa^{-1}+\kappa^{d/2}
\gamma^{1/(2\nu)})\|u_{xx}\|_{\cL_{p}},
$$
where $\nu=\mu/(\mu-1)$, $N_{1}=N(d,\delta,p)$,
 and $\kappa\geq4$ is an arbitrary number.
After   choosing 
$R=R (d,\delta,p)\in(0,R_{0}]$
and $\kappa=\kappa (d,\delta,p)\geq4$
so that 
$$
N_{1} \kappa^{-1}\leq1/4,\quad
N_{1} 
\kappa^{d/2}  
R\leq 1/4,
$$
 and finally
choosing $\gamma=\gamma(d,\delta,p)>0$ so that
$$
N_{1} \kappa^{d/2}
\gamma^{1/(2\nu)}\leq1/4,
$$
 we come to \eqref{7.1.3}.
The lemma is proved.

\mysection{Proof of Theorem \protect\ref{theorem 6.30.1}}

                                                 \label{section 7.2.4}

We take a $p\in(2,\infty)$ and take $\gamma$ from Lemma \ref{lemma 7.1.2}
and suppose that 
Assumptions   \ref{assumption 6.30.1} and
\ref{assumption 6.30.2} $(\gamma)$ are satisfied.
As usual, bearing in mind the method of continuity,
one sees that it suffices to prove the a priori estimate
\eqref{7.1.8}. 

Notice that
$$
\|L_{0}u-\lambda u\|_{\cL_{p}}
\leq \|L u-\lambda u\|_{\cL_{p}}
+N\| u_{x}\|_{\cL_{p}}
+K\| u\|_{\cL_{p}},
$$
where $N=N(d,K)$. 
Since we only consider large $\lambda$, this shows
that   it suffices to prove  
\eqref{7.1.8} with $L_{0}$ in place of $L$.
Therefore, below we assume that $b=c=0$.

In that case by using partitions of unity
one easily derives from Lemma \ref{lemma 7.1.2} that
for any $u\in W^{2}_{p}$
$$
\|u_{xx}\|_{\cL_{p}}\leq N(\|Lu\|_{\cL_{p}}+\| u_{x}\|_{\cL_{p}}
+\| u\|_{\cL_{p}}),
$$
where $N=N(d,\delta,p,R_{0})$. 
Using the interpolation inequality  
$$\|u_{x}\|_{\cL_{p}}
\leq\varepsilon\|u_{xx}\|_{\cL_{p}}+N(d,p)\varepsilon^{-1}
\|u \|_{\cL_{p}}, \quad\varepsilon>0,
$$
 shows that
\begin{equation}
                                                          \label{7.2.1}
\|u_{xx}\|_{\cL_{p}}\leq N(\|Lu\|_{\cL_{p}} 
+\| u\|_{\cL_{p}}).
\end{equation}
It follows that for
any $\lambda\geq0$
$$
\lambda \|u\|_{\cL_{p}}+\sqrt{\lambda}\|u_{x}\|_{\cL_{p}}+
\|u_{xx}\|_{\cL_{p}}
$$
$$
\leq N(\|Lu-\lambda u\|_{\cL_{p}} 
+(\lambda+1)\| u\|_{\cL_{p}}),
$$
 which implies that we only need to find
$\lambda_{0}(d,\delta,p,R_{0})
\geq1$ such that for $\lambda\geq\lambda_{0}$
we have  
\begin{equation}
                                                          \label{7.1.9}
\lambda \|u\|_{\cL_{p}} 
\leq N \|L u-\lambda u\|_{\cL_{p}} 
\end{equation}
with $N=N(d,\delta,p,R_{0} )$.

As is usual in such situations, we will follow an idea 
suggested by S.~Agmon.
 Consider the space
$$
\bR^{d+1}=\{z=( x,y): x\in\bR^{d}, y\in\bR\}
$$
 and the function
$$
\tilde{u}( z)=u(t,x)\zeta(y)\cos(\mu y),
$$
 where
 $\mu=\sqrt{\lambda}$ and $\zeta$ is a $C^{\infty}_{0}(\bR)$ function,
$\zeta\not\equiv0$.
Also introduce the operator
$$
\tilde{L}v(t,z)=a^{ij}(x)v_{x^{i}x^{j}}(z)+v_{yy}( z).
 $$

As is easy to see, the operator $\tilde{L}$ satisfies
Assumption \ref{assumption 6.30.2} ($\gamma'$) (relative to $\bR^{d+1}$)
with $\gamma'=N(d)\gamma$. Therefore, by reducing the
$\gamma$ taken from Lemma \ref{lemma 7.1.2} if necessary,
we may   apply 
the above results to the operator $\tilde{L}$ and
 in light of \eqref{7.2.1} applied to
$\tilde{u}$ and $\tilde{L}$ we get
\begin{equation}
                                                     \label{7.2.2}
\|\tilde{u}_{zz}\|_{\cL_{p}(\bR^{d+1})}
\leq N(\|\tilde{L}\tilde{u}  \|_{\cL_{p}(\bR^{d+1})}
+\|\tilde{u}\|_{\cL_{p}(\bR^{d+1})}).
\end{equation}

It is not hard to see that
$$
\int_{\bR}|\zeta(y)\cos(\mu y)|^{p}\,dy
 $$
is bounded   away from zero for $\mu\in\bR$.
Therefore,
$$
\|u\|^{p}_{\cL_{p}(\bR^{d })}
 =\mu^{-2p}\big(\int_{\bR}|\zeta(y)\cos(\mu y)|^{p}\,dy\big)^{-1}
\int_{\bR^{d+1}}\big|\tilde{u}_{yy}( z)
$$
$$
 -u( x)[\zeta''(y)\cos(\mu y)-2\mu\zeta'(y)\sin(\mu y)
 ]\big|^{p}\,dz 
$$
$$
 \leq N\mu^{-2p}\big(\|\tilde{u}_{zz}\|^{p}_{\cL_{p}(\bR^{d+1})}
+(\mu^{p}+1)\|u \|^{p}_{\cL_{p}(\bR^{d }) }\big).
$$
This and \eqref{7.2.2} yield 
$$
\mu^{2} \|u \|_{\cL_{p} }\leq 
N\|\tilde{L}\tilde{u} \|_{\cL_{p}(\bR^{d+1})}
+N(\mu+1)\|u \|_{\cL_{p}  }.
$$
Since
$$
\tilde{L}\tilde{u} =\zeta \cos(\mu y)
[Lu -\lambda u]+u[\zeta''\cos(\mu y)-2\mu\zeta'\sin(\mu y)],
$$
we have
$$
\|\tilde{L}\tilde{u}  \|_{\cL_{p}(\bR^{d+1})}\leq
N\|Lu -\lambda u\|_{\cL_{p} }
+  N(\mu+1)\|u\|_{\cL_{p}  },
$$
so that
$$
\lambda \|u \|_{\cL_{p}  }\leq 
N_{1}\|Lu -\lambda u\|_{\cL_{p} }
+  N_{2}(\sqrt{\lambda}+1)\|u\|_{\cL_{p}  }.
$$
For $\lambda\geq\lambda_{0}=16N_{2}^{2}+4N_{2}$ we have
$$
N_{2}\sqrt{\lambda}\leq(1/4)\lambda,\quad
N_{2}\leq(1/4)\lambda,\quad
 N_{2}(\sqrt{\lambda}+1)\leq(1/2)\lambda
$$
and we arrive at \eqref{7.1.9} with $N=2N_{1}$.
The theorem is proved.

\end{document}